\documentclass[12pt]{article}
\usepackage{amsmath,amssymb,amsfonts,amsthm}
\usepackage{epsfig}

\newtheorem{theorem}{Theorem}[section]

\newtheorem{proposition}[theorem]{Proposition}

\newtheorem{corollary}[theorem]{Corollary}
\newtheorem{conjecture}{Conjecture}

\newtheorem*{Main}{Main Theorem}

\def\pfsp{\hskip 1em}

\newcommand{\skpf}[1]{\noindent{\it Sketch of Proof\,:}\pfsp #1 {\hfill $\Box$} \smallskip}

\def \QD1 {\hfill $\spadesuit$}

\newcommand{\DF}[1]{{\bf #1\/}}

\newcommand{\set}[2]{\{#1 \;|\; #2 \}}
\newcommand{\ems}{\varnothing}
\newcommand{\sm}{\setminus}

\newcommand{\cn}{\chi}

\newcommand{\cGD}{{\cal DG}}

\newcommand{\cD}{{\cal D}}

\newcommand{\CR}{{\mathrm{CRI}}}
\newcommand{\vCR}{\overrightarrow{\mathrm{CRI}}}
\newcommand{\EX}{{\mathrm{EXT}}}
\newcommand{\vEX}{\overrightarrow{\mathrm{EXT}}}
\newcommand{\ex}{\mathrm{ext}}
\newcommand{\vex}{\overrightarrow{\mathrm{ext}}}
\newcommand{\vDG}{\overrightarrow{\cal DG}}

\newcommand{\dicn}{\vec{\chi}}

\newcommand{\f}{\varphi}
\newcommand{\fin}{\varphi^{-1}}

\newcommand{\nat}{\mathbb{N}}

\newcommand{\ganz}{\mathbb{Z}}

\newcommand{\dis}{\boxplus}

\newcommand{\mdo}[3]{#1 \equiv #2 \, \mathrm{(mod} \; #3 \mathrm{)}}

\newcommand{\NP}  {{\sf NP}}
\newcommand{\NPC} {{\NP}-complete}

\newcommand{\jou}[4]{{\rm #1} #2 (#3) #4.}

\def \DM {Discrete Math.}
\def \SIDMA {SIAM J. Discrete Math.}
\def \SIC {SIAM J. Comput.}

\numberwithin{equation}{section}

\begin{document}
\title{\bf Minimum number of arcs in $k$-critical digraphs with order at most $2k-1$}

\author{Lucas Picasarri-Arrieta\thanks{Universit\'e C\^ote d'Azur, CNRS, Inria, I3S, Sophia-Antipolis, France (email:lucas.picasarri-arrieta@inria.fr), financial supports: DIGRAPHS ANR-19-CE48-0013 and EUR DS4H Investments ANR-17-EURE-0004.} \and Michael Stiebitz\thanks{Technische Universit\"at Ilmenau. Inst. of Math. PF 100565, D-98684 Ilmenau, Germany (email:michael.stiebitz@tu-ilmenau.de)}
}

\date{}
\maketitle
\begin{abstract}
The dichromatic number $\dicn(D)$ of a digraph $D$ is the least integer $k$ for which $D$ has a coloring with $k$ colors such that there is no monochromatic directed cycle in $D$. A digraph $D$ is called $k$-critical if each proper subdigraph $D'$ of $D$ satisfies $\dicn(D')<\dicn(D)=k$. For integers $k$ and $n$, let $\vex(k,n)$ denote the minimum number of arcs possible in a $k$-critical digraph of order $n$. It is easy to show that $\vex(2,n)=n$ for all $n\geq 2$, and $\vex(3,n)\geq 2n$ for all possible $n$, where equality holds if and only if $n$ is odd and $n\geq 3$. As a main result we prove that if $n, k$ and $p$ are integers with $n=k+p$ and $2\leq p \leq k-1$, then $\vex(k,n)=2({\binom{n}{2}} - (p^2+1))$, and we give an exact characterisation of $k$-critical digraphs for which equality holds.
This generalizes a result about critical graphs obtained in 1963 by Tibor Gallai.
\end{abstract}



\section{Introduction and main result}

This paper is concerned with the coloring problem for digraphs.
Let $D$ be a digraph. An \DF{acyclic $k$-coloring} of $D$ is a map $\f:V(D)\to C$ with color set $C=\{1,2, \ldots, k\}$ such that for each color $c\in C$ the subdigraph of $D$ induced by the color class $\fin(c)=\set{v\in V(D)}{\f(v)=c}$ is acyclic. The \DF{dichromatic number} of $D$, denoted by $\dicn(D)$, is the least integer $k\geq 0$ for which $D$ admits an acyclic $k$-coloring. This coloring concept for digraphs, introduced in 1982 by
Neumann-Lara and Erd\H{o}s (cf. \cite{ErdosNL82,NeumannLara82}), generalizes the ordinary coloring concepts for graphs and has been studied in numerous papers until today (cf. \cite{BangJensenBSS2019,BokalFJKM2004,Erdos79,Golowich2016,HarayunyanMo2011,LiMohar2016,Mohar2010,NeumannLara85,NeumannLara94}).

Let $k$ be a non-negative integer. We call a digraph $D$ \DF{critical} and \DF{$k$-critical} if every proper subdigraph $D'$ of $D$ satisfies $\dicn(D')<\dicn(D)=k$. Note that any digraph $D$ has a critical subdigraph with the same dichromatic number as $D$. Hence several problems related to the dichromatic number can be reduced to critical digraphs. Critical digraphs were first introduced and investigated by Neumann-Lara \cite{NeumannLara82}. Let $\vCR(k)$ denote the class of $k$-critical digraphs, and for an integer $n$, let
$$\vCR(k,n)=\set{D\in \vCR(k)}{|D|=n
}.$$

Note that $\dicn(D)=0$ if and only if $D=\ems$, and $\dicn(D)\leq 1$ if and only if $D$ is acyclic.
Consequently, $\vCR(0)=\{\ems\}$, $\vCR(1)=\{K_1\}$, and $\vCR(2)$ is the class of directed cycles. In 2007, Chen, Hu, and Zang~\cite{ChenHZ2007} proved that the decision problem whether a given digraph $D$ satisfies $\dicn(D)\leq 2$ is {\NPC}. Hence a characterization of the class $\vCR(k)$ for fixed $k\geq 3$ is unlikely. In this paper, we are interested in the extremal function $\vex(k,n)$ defined by
$$\vex(k,n)=\min \set{|A(D)|}{D\in \vCR(k,n)}$$
as well as the corresponding class
$$\vEX(k,n)=\set{D\in \vCR(k,n)}{|A(D)|=\vex(k,n)}$$
of extremal graphs. The following theorem, which is the main result of this paper, implies a classical result for critical graphs that was obtained by Gallai \cite{Gallai63b} in the 1960s.

\begin{Main}
Let $n$, $k$ and $p$ be integers satisfying $n=k+p$ and $2\leq p \leq k-1$. Then $\vex(k,n)=2({\binom{n}{2}} - (p^2+1))$.
\end{Main}

We also give an exact characterization of $\vEX(k,n)$ under the same assumptions on $k,n$.

\section{Preliminaries}

Our notation is standard. For integers $k$ and $\ell$, let $[k,\ell]=\set{x\in \ganz}{k \leq x \leq \ell}$.

\subsection{Notation for graphs and digraphs}

The term \DF{graph} in this paper always refers to a finite undirected graph without multiple edges and without loops.  For a graph $G$, we denote by $V(G)$ and $E(G)$ the \DF{vertex set} and the \DF{edge set} of $G$, respectively. The \DF{order} of $G$, denoted by $|G|$, is the number of vertices of $G$. If $|G|=0$, then $G$ is called the \DF{empty graph}, briefly $G=\ems$. As usual, $K_n$ denotes the complete graph of order $n$, and $C_n$ denotes the cycle of order $n$ with $n\geq 3$. A cycle is called \DF{odd} or \DF{even} depending on whether its order is odd or even.
A \DF{proper $k$-coloring} of $G$ is a map $\f : V(G) \to [1,k]$ such that for each color $c\in [1,k]$, the color class $\f^{-1}(c)$ induces an independent set. The \DF{chromatic number} of $G$, denoted by $\chi(G)$, is the least integer $k$ for which $G$ admits a proper $k$-coloring. 

Digraphs are always assumed to be finite and simple, that is, no loops and no parallel arcs in the same direction are allowed. However, digraphs may have \DF{digons}, that is, two arcs with opposite directions between the same pair of vertices. For a digraph $D$, we denote the \DF{vertex set} of $D$ by $V(D)$, and the \DF{arc set} of $D$ by $A(D)$. The number of vertices of $D$ is called the \DF{order} of $D$, denoted by $|D|$. A digraph $D$ with $|D|=0$ is said to be \DF{empty}, in this case we also write $D=\ems$. A digraph is said to be \DF{acyclic} if it contains no directed cycle as a subdigraph.

Let $u$ and $v$ be two distinct vertices of a digraph $D$.
Let $uv$ be an arc of $D$. We say that $u$ is a \DF{neighbor} of $v$ and vice versa. We also say that $v$ is an \DF{out-neighbor} of $u$ and $u$ is an \DF{in-neighbor} of $v$. By $N_D^+(v)$ we denote the set of out-neighbors of $v$; by $N_D^-(v)$ the set of in-neighbors of $v$. Then $d_D^+(v)=|N_D^+(v)|$ is the \DF{out-degree} of $v$ in $D$, and $d_D^-(v)=|N_D^-(v)|$ is the \DF{in-degree} of $v$ in $D$.

We denote by $D[X]$ the subdigraph of $D$ that is \DF{induced} by the vertex set $X$, that is, $V(D[X])=X$ and $A(D[X])= \set{uv \in A(D)}{u,v\in X}$. A digraph $D'$ is called an \DF{induced subdigraph} of $D$ if $D'=D[V(D')]$.  If $X$ is a subset of $V(D)$, then we define $D-X=D[V(D) \setminus X]$. If $X=\{v\}$ is a singleton, we use $D-v$ rather than $D-X$. For an arc set $F\subseteq A(D)$, let $D-F$ denote the subdigraph of $D$ with vertex set $V(D)$ and arc set $A(D)\sm F$. If $F=\{a\}$ is a singleton, we write $D-a$ instead of $D-F$.
The \DF{directed cycle} of order $n\geq 2$, denoted by $\vec{C}_n$, is the digraph with vertex set $\{v_1,\dots,v_n\}$ and arc set $\{v_iv_{i+1}\mid i\in [1,n]\}$ where subscripts are taken modulo $n$. The digraph $\vec{C}_n$ is also called the \DF{$n$-cycle}. A digraph $D$ is called \DF{acyclic} if $D$ contains no directed cycle as a subdigraph. 

The \DF{bidirected graph} associated to a graph $G$, denoted by $D^\pm(G)$, is the digraph obtained from $G$ by replacing each edge of $G$ by a digon between the same vertices. So $D^\pm(K_1)=K_1$, and $D^\pm(K_2)=\vec{C}_2$ is a directed cycle of order $2$. 
Observe that, if $G$ is a graph and $D= D^\pm(G)$, then there is a one-to-one correspondence between the proper colorings of $G$ and the acyclic colorings of $D$ and, in particular, $\chi(G)=\dicn(D^\pm(D))$.

\subsection{Minimum number of edges in critical graphs}

Deleting a vertex or edge of a graph decreases its chromatic number by at most one. Let $k$ be a non-negative integer. A graph $G$ is called \DF{critical} and \DF{$k$-critical} if every proper subgraph $G'$ of $G$ satisfies $\cn(G')<\cn(G)=k$. Let $\CR(k)$ denote the class of $k$-critical graphs, and for $n\geq 0$, let
$$\CR(k,n)=\set{G\in \CR(k)}{|G|=n}.$$
Clearly, if $G\in \CR(k)$, then $|G|\geq k$ and equality holds if and only if $G=K_k$. For $k\leq 2$, $K_k$ is the only $k$-critical graph. In particular, $\CR(0)=\{\ems\}$, $\CR(1)=\{K_1\}$, and $\CR(2)=\{K_2\}$. Furthermore, K\"onig's theorem, saying that a graph $G$ is bipartite if and only if $G$ contains no odd cycle as a subgraph, implies that
$$\CR(3)=\set{C_n}{\mdo{n}{1}{2}}.$$
A good characterization of $\CR(k)$ with fixed $k\geq 4$ is very unlikely. It is known that if $k\geq 4$, then $\CR(k,n)\not=\ems$ if and only if $n\geq k$ and $n\not= k+1$. The concept of critical graphs is due to G. A. Dirac. In the 1950s he established the basic properties of critical graphs in a series of papers (see \cite{Dirac52,Dirac53,Dirac57}), and the study of critical graphs has attracted a lot of attention since then. Dirac also started to investigate the function $\ex(k,n)$ defined by
$$\ex(k,n)=\min \set{|E(G)|}{G\in \CR(k,n)}$$
as well as the corresponding class
$$\EX(k,n)=\set{G\in \CR(k,n)}{|E(G)|=\ex(k,n)}$$
of extremal graphs.
For a detailed discussion of the many partial results related to this extremal function the reader is referred to the survey by Kostochka \cite{Kostochka2006} (see also \cite{KostochkaY14a,Stiebitz2024}).

It is easy to show that every graph in $\CR(k)$ with $k\geq 1$  has minimum degree at least $k-1$, which gives the trivial lower bound $2\ex(k,n)\geq (k-1)n$. By Brooks' famous theorem \cite{Brooks41} it follows that $2\ex(k,n)=(k-1)n$ if and only if $n=k$ or $k=3$ and $n$ is odd. In 1957 Dirac \cite{Dirac57} (see also \cite{Dirac74}) proved that $2\ex(k,n)\geq (k-1)n+k-3$ for every $k\geq 4$ and $n\geq k+2$. For $k\geq 3$, we denote by $\cGD(k)$ the class of graphs that can be obtained from two disjoints $K_k$'s by splitting a vertex of one $K_k$ into an edge of the other $K_k$. So a graph $G$ belongs to $\cGD(k)$ if and only if the vertex set of $G$ consists of three non-empty pairwise disjoint sets $X, Y_1$ and $Y_2$ with
    $$|Y_1|+|Y_2|=|X|+1=k-1$$
and two additional vertices $v_1$ and $v_2$ such that $G[X]$ and $G[Y_1 \cup Y_2]$ are complete graphs not joined by any edge in $G$, and $X \cup Y_i$ is the neighborhood of $v_i$ in $G$ for $i\in \{1,2\}$. Then it is easy to show that $\cGD(k)\subseteq \CR(k,2k-1)$. The class $\cGD(k)$ was introduced and investigated by Dirac \cite{Dirac74} and by Gallai \cite{Gallai63a}. Dirac \cite{Dirac74} proved that $\EX(k,2k-1)=\cGD(k)$ for every $k\geq 3$. In 1963, Gallai published two fundamental papers \cite{Gallai63a,Gallai63b} about the structure of critical graphs. In particular, he proved the following two remarkable results.

\begin{theorem} [Gallai]
\label{theorem:gallai1}
Let $G$ be a $k$-critical graph of order $n$. If $n\leq 2k-2$, then the complement of $G$ is disconnected.
\end{theorem}

\begin{theorem} [Gallai]
\label{theorem:gallai2}
Let $n=k+p$ be an integer, where $k,p\in \nat$ and $2\leq p \leq k-1$.
Then $$\ex(k,n)={\binom{n}{2}} - (p^2+1) $$
and $\EX(k,n)$ is the class of graphs $G$ that can be obtained from two disjoint subgraphs $G_1$ and $G_2$ with $G_1=K_{k-p-1}$ and $G_2\in \cGD(p+1)$ by joining each vertex of $G_1$ to each vertex of $G_2$ by exactly one edge.
\end{theorem}

In 2014 Kostochka and Yancey \cite{KostochkaY14a} established a lower bound for $\ex(k,n)$ that is sharp when $k\geq 4$ and $\mdo{n}{1}{(k-1)}$. In particular, they proved that if $k\geq 4$, then
\begin{align}
\label{Equation:KostochkaYancey-limit}
\lim_{n \to \infty}\frac{\ex(k,n)}{n}=\frac{1}{2}(k-\frac{2}{k-1}).
\end{align}
Furthermore, Kostochka and Yancey \cite{KostochkaY14b} proved that

\begin{align}
\label{Equation:KostochkaYancey-4critica}
\ex(4,n)=\left\lceil \frac{5n-2}{3}\right\rceil
\end{align}
As pointed out in \cite{KostochkaY14b} the above result and its short proof leads to a short proof of Gr\"otzsch's famous Dreifarbensatz \cite{Grotzsch58} saying that any planar triangle free graph $G$ satisfies $\cn(G)\leq 3$.

\subsection{Minimum number of arcs in critical digraphs}

While the study of critical graphs has attracted a lot of attention, not so much is known about the structure of critical digraphs. 
One can easily show that, in a $k$-critical digraph $D$, every vertex $v$ satisfies $d^+_D(v) \geq k-1$ and $d^-_D(v) \geq k-1$.
Mohar~\cite{Mohar2010} established a Brooks-type result for digraphs, characterizing exactly the $k$-critical digraphs in which every vertex has in- and out-degree exactly $k-1$. 
From this result we deduce that, for every $k,n \in \nat$, we have $\vex(k,n)\geq (k-1)n$ and equality holds if and only if $k=2$ and $n\geq 2$, or $k=3$ and $n\geq 3$ is odd, or $n=k$. Furthermore, $\vEX(2,n)=\{\vec{C}_n\}$ for $n\geq 2$, $\vEX(3,n)=\{D^\pm(C_n)\}$ for odd $n$ and $n\geq 3$, and $\vEX(k,k)=\{D^\pm(K_k)\}$.

Let $D=D^\pm(G)$ be the bidirected graph associated to a graph $G$. Then $\cn(G)=\dicn(D)$ and it is easy to see that $G\in \CR(k)$ if and only if $D\in \vCR(k)$. As a consequence, we obtain $\vex(k,n)\leq 2\ex(k,n)$ if $n\geq k \geq 4$ and $n\not=k+1$. Note that $\CR(k,k+1)=\ems$, but $\vCR(k,k+1)\not=\ems$ (see Proposition~\ref{prop:case_p_1}). Kostochka and Stiebitz \cite{KostochkaS2019} proposed the following conjecture.

\begin{conjecture} [Kostochka and Stiebitz]
\label{conjecture:kost+stieb}
Let $k,n\in \nat$ with $n\geq k\geq 4$ and $n\not=k+1$. Then $\vex(k,n)=2\ex(k,n)$ and hence
$$\lim_{n \to \infty}\frac{\vex(k,n)}{n}=(k-\frac{2}{k-1}).$$
Furthermore, $\vEX(k,n)=\set{D^\pm(G)}{G\in \EX(k,n)}$.
\end{conjecture}

A first result supporting this conjecture was obtained by Kostochka and Stiebitz \cite{KostochkaS2019}; they proved that if $n\geq 4$, then
$$\vex(4,n)\geq \frac{10n-4}{3}.$$

Havet, Picasarri-Arrieta, and Rambaud~\cite{Havet2023} later characterized the case of equality. They proved that if $D\in \vCR(4,n)$ has $\frac{10n-4}{3}$ arcs, then $D \in \set{D^\pm(G)}{G\in \EX(4,n)}$.
The best known result towards Conjecture~\ref{conjecture:kost+stieb} is due to Aboulker and Vermande (see Theorem~10 in~\cite{Vermande2022}), who showed that
\[
(k-\frac{1}{2}-\frac{1}{k-1}) \leq \lim_{n \to \infty}\frac{\vex(k,n)}{n} \leq (k-\frac{2}{k-1}).
\]

In the proof of the main theorem we will use two recent results. 
For each $k\geq 2$, let us define $\vDG(k)$ as follows: $\vDG(2) = \{\vec{C_3}\}$ and $\vDG(k) = \{ D^{\pm}(G) \mid G\in \cGD(k)\}$ for every $k\geq 3$.
The first result we need, due to Aboulker and Vermande (see Theorem~4 in \cite{Vermande2022}), is the following generalization of Dirac's result to digraphs.

\begin{theorem}[Aboulker and Vermande]
	\label{theorem:Dirac_digraphs}
	If $D\in \vCR(k,n)$ and $n>k\geq 4$, then
	\[
		|A(D)| \geq (k-1)n + k-3.
	\]
	Moreover, equality holds if and only if $D\in \vDG(k)$.
\end{theorem}

Observe that the result above implies our main theorem when $p=k-1$.
The second result we need is due to Stehl\'ik \cite{Stehlik2019} and deals with decomposable critical digraphs. We introduce it in the next section.

\subsection{Decomposable critical digraphs}

As pointed out by Bang-Jensen, Bellitto, Schweser, and Stiebitz \cite{BangJensenBSS2019}, two well known constructions for critical graphs have counterparts for digraphs, namely the Haj\'os join and the Dirac join.

Let $D_1$ and $D_2$ be two disjoint digraphs. Let $D'$ be the digraph obtained from the union $D_1 \cup D_2$ by adding all possible arcs in both directions between $D_1$ and $D_2$, that is, $V(D')=V(D_1) \cup V(D_2)$ and $A(D)=A(D_1) \cup A(D_2) \cup \{uv, vu ~|~ u \in V(D_1) \text{ and } v \in V(D_2)\}$. We call $D'$ the \DF{Dirac join} of $D_1$ and $D_2$ and write $D'= D_1 \dis D_2$. A proof of the following is straightforward.

\begin{theorem}[Dirac Construction]
Let $D=D_1 \dis D_2$ be the Dirac join of two disjoint non-empty digraphs $D_1$ and $D_2$. Then, $\chi(D) = \chi(D_1) + \chi(D_2)$ and $D$ is critical if and only if both $D_1$ and $D_2$ are critical.
\label{theorem:dirac-join}
\end{theorem}

A digraph is called \DF{decomposable} if it is the Dirac join of two non-empty subdigraphs; otherwise the digraph is called \DF{indecomposable}. Theorem~\ref{theorem:dirac-join} implies that a decomposable critical digraph is the Dirac join of its indecomposable critical subdigraphs. In this sense, the indecomposable critical digraphs are the building elements of critical digraphs. In 2019, Stehl\'ik \cite{Stehlik2019} proved the following result, thereby answering a question proposed in \cite{BangJensenBSS2019}. Note that this result applied to bidirected graphs implies the decomposition result for critical graphs due to Gallai (see Theorem~\ref{theorem:gallai1}). Stehl\'ik's proof uses matching theory. However, one can also use the hypergraph version of Theorem~\ref{theorem:gallai1} obtained by Stiebitz, Storch, and Toft \cite{StiebitzST2016} (see also~\cite[Theorem~7.25]{Stiebitz2024}).

\begin{theorem}[Stehl\'ik]
If $D$ is an indecomposable critical digraph, then $|D|\geq 2\dicn(D)-1.$
\label{theorem:decomposable}
\end{theorem}
\skpf{Let $D\in \vCR(k)$ with $|D|\leq 2k-2$. Our aim is to show that $D$ is decomposable. Let $H$ be a hypergraph with $V(H)=V(D)$ such that $e\subseteq V(H)$ is an edge of $H$ if and only if $D[e]$ is a directed cycle. Note that no edge of $H$ contains another edge of $H$. Then $\cn(H)=k$ and $\cn(H-v)<k$ for every vertex $v$ of $H$. Then $H$ has a subhypergraph $H'$ such that $\cn(H'')<\cn(H')=\cn(H)=k$ for every proper subhypergraph $H''$ of $H'$. Furthermore, we have that $V(H')=V(H)$. Then, by \cite[Main Theorem 1]{StiebitzST2016}, $H'$ is obtained from the disjoint union of two nonempty subhypergraphs by adding all $2$-edges between these two subhypergraphs. Since $e$ is a 2-edge of $H$ if and only if $D[e] = D^\pm(K_2)$, we obtain that $D$ is decomposable.}

The complement of a digraph $D$, denoted by $\overline{D}$, is the digraph with $V(\overline{D})=V(D)$ and $A(\overline{D})=\set{uv}{uv \not\in A(D)}$. So if $D$ has order $n$, then $D\cup \overline{D}=D^\pm(K_n)$. Clearly, a non-empty digraph $D$ is indecomposable if and only if $\overline{D}$ is connected.

\section{Critical digraphs whose order is near to $\dicn$}

Let $D$ be a digraph. A non-empty subdigraph $D'$ of $D$ is called a \DF{dominating subdigraph} of $D$ if there exists a non-empty subdigraph $D''$ of $D$ such that $D=D'\dis D''$. A vertex $v$ of $D$ is called a \DF{dominating vertex} of $D$ if the subdigraph $D[\{v\}]$ is a dominating subdigraph of $D$. Note that if $X\subseteq V(D)$ is a set of $p\geq 0$ vertices such that any vertex of $X$ is a dominating vertex of $D$, then $D=D^\pm(K_p)\dis (D-X)$.

\begin{theorem}
Let $D$ be a $k$-critical digraph with $k\geq 1$, let $p$ be the number of dominating vertices of $D$, and $q$ be the number of dominating directed cycles of $D$ having order at least three. Then the following statements hold:
\begin{itemize}
\item[{\rm (a)}] $0\leq p \leq k$ and there exists a digraph $D'\in \vCR(k-p)$ such that $D=D^\pm(K_p)\dis D'$, $D'$ has no dominating vertex, and $|D'| \geq \tfrac{3}{2}(k-p).$ Furthermore, $p\geq 3k-2|D|$ and equality holds if and only if $D'$ is the Dirac join of $\tfrac{1}{2}(k-p)$ disjoint directed $3$-cycles of $D$.
\item[{\rm (b)}] $0 \leq p+2q\leq k$ and there exist digraphs $D_1\in \vCR(2q)$ and $D_2\in \vCR(k-p-2q)$
such that $D=D^\pm(K_p) \dis D_1 \dis D_2$, $D_1$ is the Dirac join of $q$ directed cycles of $D$ each of which has order at least three, $D_2$ has no dominating vertex and no dominating directed cycle, and $|D_2|\geq \tfrac{5}{3}(k-p-2q)$. Furthermore, $2p+q\geq 5k-3|D|$ and equality holds if and only if
$D_1$ is the Dirac join of $q$ directed $3$-cycles of $D$, and $D_2$ is the Dirac join of $\tfrac{1}{3}(k-p-2q)$ disjoint subdigraphs of $D$ each of which belong to $\vCR(3,5)$.
\end{itemize}
\label{theorem:n+near+to:dicn}
\end{theorem}
\begin{proof}
Note that $0 \leq p+2q \leq k$ trivially holds as the Dirac join of $D^\pm(K_{p})$ and $q$ directed cycles has dichromatic number $p+2q$ by Theorem~\ref{theorem:dirac-join}.
In what follows, let $D\in \vCR(k)$ with $k\geq 1$ and let
$$D=D_1 \dis D_2 \dis \cdots  \dis D_s,$$
where $\overline{D}_1, \overline{D}_2, \ldots, \overline{D}_s$ are the connected components of $\overline{D}$. Obviously, $s\geq 1$. For $i\in [1,s]$, let $k_i=\dicn(D_i)$ and $n_i=|D_i|$. By Theorem~\ref{theorem:dirac-join}, we obtain that
\begin{itemize}
\item[{\rm (1)}] $k=k_1+k_2 + \cdots + k_s$ and $D_i\in \vCR(k_i,n_i)$ for $i\in [1,s]$.
\end{itemize}
Since $\overline{D}_i$ is connected, Theorem~\ref{theorem:decomposable} implies that
\begin{itemize}
\item[{\rm (2)}] $n_i \geq 2k_i-1$ for $i\in [1,s]$.
\end{itemize}
Since $\vCR(1)=\{K_1\}$, and $\vCR(2)=\set{\vec{C}_n}{n\geq 2}$, we obtain that $k_i=1$ and $D_i=K_1$, or $k_i=2$ and $n_i\geq 3$ (where equality holds if and only if $D_i=\vec{C}_3$), or $k_i\geq 3$ and $n_i\geq 5$. For a subset $I$ of $[1,s]$, let $D_I=\dis_{i\in I}D_i$ be the Dirac join of the digraphs $D_i$ with $i\in I$, and let $k_I=\sum_{i\in I}k_i,$ where $D_{\ems}=\ems$ and $k_{\ems}=0$.
Then it follows from Theorem~\ref{theorem:dirac-join} that $D_I\in \vCR(k_I)$. Let $P=\set{i\in [1,s]}{k_i=1}$, $Q=\set{i\in [1,s]}{k_i=2}$, $R=[1,s]\sm (P \cup Q)$, $p=|P|$, $q=|Q|$, and $r=|R|$. Then $P,Q$ and $R$ are pairwise disjoint sets whose union is $[1,s]$. Thus we obtain that
\begin{itemize}
\item[{\rm (3)}] $D=D_P\dis D_Q \dis D_R$, where $D_P=D^\pm(K_p)$ and $D_Q\in \vCR(2q)$.
\end{itemize}
Note that $p$ is the number of dominating vertices of $D$, and $q$ is the number of dominating directed cycles of $D$.

To establish a lower bound for $p$, let $\overline{P}=[1,s] \sm P$. Then $\overline{P}=R \cup Q$ and $D=D_P\dis D_{\overline{P}}$.
For $i\in \overline{P}$, we have that $k_i\geq 2$ and so, by (2), $n_i\geq 2k_i-1\geq \tfrac{3}{2}k_i$, where equality holds if and only if $D_i=\vec{C}_3$. From Theorem~\ref{theorem:dirac-join} and (1) it follows that $k_P=p$ and $k_{ \overline{P}}=k-p$. For the order of $D$, it then follows from (1) and (2) that
$$|D|=p+\sum_{i\in \overline{P}}n_i\geq p+\tfrac{3}{2}\sum_{i\in \overline{P}}k_i= p+\tfrac{3}{2}(k-p),$$
which is equivalent to $p\geq 3k-2|D|$. Clearly, $p= 3k-2|D|$ if and only if $D_{\overline{P}}$ is the Dirac join of $\tfrac{1}{2}(k-p)$ disjoint $\vec{C}_3$'s. This proves (a).

\medskip

For $i\in R$, we have $k_i\geq 3$ and so, by (2), $n_i\geq 2k_i-1\geq \tfrac{5}{3}k_i$, where equality holds if and only if $D_i\in \vCR(3,5)$.
From Theorem~\ref{theorem:dirac-join} we conclude that $k_P=p$, $k_Q=2q$, and $k_R=k-p-2q$.
For the order of $D$ we then obtain that
$$|D|=p+\sum_{i\in Q}n_i+\sum_{i\in R}n_i\geq p+3q+\tfrac{5}{3}\sum_{i\in R}k_i=p+3q+\tfrac{5}{3}(k-p-2q),$$
which is equivalent to $2p+q\geq 5k-3|D|$. Clearly, $2p+q= 5k-3|D|$ if and only if $D_i=\vec{C}_3$ for all $i\in Q$ and $D_i\in \vCR(3,5)$ for all $i\in R$. Thus (b) is proved.
\end{proof}

\medskip

In the remaining of this section, we discuss some consequences of Theorem~\ref{theorem:n+near+to:dicn}.

For a digraph $K$ and a class of digraphs $\cD$, define $K\dis \cD=\set{K\dis D}{D\in \cD}$. If $\cD$ is a digraph property (i.e., a class of digraphs closed under taking isomorphic copies), then we do not distinguish between isomorphic digraphs, so we are only interested in the number of isomorphism types of $\cD$, that is, the number of equivalence classes of $\cD$ with respect to the isomorphism relation for digraphs. Clearly, the number of isomorphisms types of the class $\vCR(k,n)$ is finite.

\medskip

From Theorem~\ref{theorem:n+near+to:dicn}(a) it follows that $\vCR(k,k+1)=D^\pm(K_{k-2})\dis \vCR(2,3)$ if $k\geq 3$. Since $\vCR(2,3)=\{\vec{C}_3\}$, we have the following result (see Lemma 16 of~\cite{Vermande2022} for an alternative proof).

\begin{proposition}
	\label{prop:case_p_1}
	For every integer $k\geq 2$,
    \[
    		\vCR(k,k+1)=\{D^\pm(K_{k-2})\dis \vec{C}_3\}.
    \]

\end{proposition}

For $k\geq 4$, Theorem~\ref{theorem:n+near+to:dicn}(a) implies that a $k$-critical digraph on $k+2$ vertices contains at least $k-4$ dominating vertices. Furthermore, if it contains exactly $k-4$ dominating vertices, Theorem~\ref{theorem:n+near+to:dicn}(b) implies that it also contains two disjoint dominating directed triangles. We thus obtain the following result.
\begin{corollary}
	For every integer $k\geq 4$,
    \[
    		\vCR(k,k+2) = \Big(D^\pm(K_{k-4})\dis \vec{C}_3 \dis \vec{C}_3\Big) \cup \Big(D^\pm(K_{k-3})\dis \vCR(3,5)\Big).
    \]
\end{corollary}

Let $\vCR^*(k,n)$ be the class consisting of all digraphs in $\vCR(k,n)$ with no dominating vertex.
For $k$-critical digraphs on $k+p$ vertices, the following is a nice consequence of Theorem~\ref{theorem:n+near+to:dicn}. In particular, it implies that the number of $k$-critical digraphs on $k+p$ vertices (up to isomorphism) is bounded by a function depending only on $p$.

\begin{corollary}
	\label{cor:caracterisation_k_p}
	For every integer $k,p$ such that $k > 2p \geq 4$,
    \[
    		\vCR(k,k+p) = \bigcup_{\ell=2}^{2p} \left\{ D^{\pm}(K_{k-\ell}) \dis  \vCR^*(\ell,\ell+p) \right\}
    \]
\end{corollary}
\begin{proof}
	Let $D$ be a $k$-critical on $k+p$ vertices. Then $|D| = k+p < \frac{3}{2}k$. Let $s$ be the number of dominating vertices of $D$ and $\ell = k-s$, so $D=D^{\pm}(K_s) \dis D'$ where $D'$ belongs to $\vCR^*(\ell, \ell+p)$. It remains to show that $\ell \in [2,2p]$.
	
	Since $p\geq 2$, we obviously have $\ell \geq 2$. On the other hand, $D'$ contains no dominating vertex, so $|D'| \geq \frac{3}{2}\ell$ (by applying Theorem~\ref{theorem:n+near+to:dicn}(a) on $D'$). We deduce
	\[
		|D'| = \ell+p \geq \frac{3}{2}\ell.
	\]
	This implies $\ell \leq 2p$ as desired.
\end{proof}

\section{Proof of the main theorem}

The following result clearly implies the Main Theorem.

\begin{theorem}
Let $p,k,n\in \nat$ such that $1\leq p \leq k-1$ and $n=k+p$. Then
\[
\vEX(k,n) = D^{\pm}(K_{k-p-1}) \dis \vDG(p+1).
\]
As a consequence, we have
\[
\vex(k,n) = \left\{
         \begin{array}{ll}
             2{\binom{n}{2}} - 3 & \mbox{if } p = 1\\
             2\left({\binom{n}{2}} - (p^2+1)\right) & \mbox{otherwise. }\\
		\end{array}
     \right.
\]

\label{theorem:Gallai+digraph}
\end{theorem}

\begin{proof}
We proceed by induction on $k$. If $k=2$, then $p=1$ and we have $\vEX(2,3)=\{\vec{C}_3\}=\vDG(2)$. If $k=3$, then $p\in \{1,2\}$ and we have $\vCR(3,4)=\{D^\pm(K_1)\dis \vec{C}_3\}=\vEX(3,4)$ (by Proposition~\ref{prop:case_p_1}) and $\vEX(3,5)=\vDG(3)$ (by Theorem~\ref{theorem:Dirac_digraphs}). So Theorem~\ref{theorem:Gallai+digraph} holds for $k=2$ and $k=3$.

We now assume $k\geq 4$. We also assume that $2 \leq p \leq k-2$, for otherwise $p\in \{1,k-1\}$ and the result is a consequence of Proposition~\ref{prop:case_p_1} or Theorem~\ref{theorem:Dirac_digraphs}. Let $n=k+p$ and let $D\in \vEX(k,n)$.
We will show that $D$ belongs to the digraph class $D^{\pm}(K_{k-p-1}) \dis \vDG(p+1)$.
	
	\medskip
	
	\textbf{Case 1:} \textit{$D$ contains a dominating vertex $v$.} Then $D=K_1\dis D'$ and $D'=D-v\in \vCR(k-1,n-1)$ (by Theorem~\ref{theorem:dirac-join}). Furthermore, $D'\in \vEX(k-1,n-1)$, for otherwise there is a digraph $\tilde{D}\in \vEX(k-1,n-1)$ with $|A(\tilde{D})| < |A(D')|$ implying that $K_1\dis \tilde{D}\in \vCR(k,n)$ has fewer  arcs than $D$, which is impossible. Since $n\leq 2k-2$, we have $(k-1)+1 \leq n-1\leq 2(k-1) - 1$, which allows us to apply the induction hypothesis to $D'$. Hence $D'$ belongs to $D^{\pm}(K_{k-p-2}) \dis \vDG(p+1)$. Consequently, $D=K_1\dis D'$ belongs to $D^{\pm}(K_{k-p-1}) \dis \vDG(p+1)$ as desired.
	
	\medskip

	\textbf{Case 2:} \textit{$D$ contains a dominating $\vec{C_3}$ but no dominating vertex.} Then $D=\vec{C_3}\dis D'$ and $D'\in \vCR(k-2,n-3)$ (by Theorem~\ref{theorem:dirac-join}). Note also that $n-3\geq k-1$, for otherwise $D'$ is the complete digraph on $k-2$ vertices, which implies that $D$ contains a dominating vertex. Hence we may apply the induction to $D'$, which implies that $D'$ contains at least $2\left({\binom{n-3}{2}}-((p-1)^2+1)\right)$ arcs. Consequently, we obtain  that
	\begin{align*}
		|A(D)| &= |A(D')| + 3 + 6(n-3)\\
			&\geq  2\left({\binom{n-3}{2}}-((p-1)^2+1)\right) + 3 + 6(n-3)\\
			&= 2\left({\binom{n}{2}} -(p^2+1)\right) -5 + 4p\\
			&\geq 2\left({\binom{n}{2}} -(p^2+1)\right) + 3,
	\end{align*}
	where in the last inequality we used $p\geq 2$.
	This is a contradiction to the minimality of $D$, since every digraph in $D^{\pm}(K_{k-p-1}) \dis \vDG(p+1)$ contains exactly $2\left( {\binom{n}{2}} - (p^2+1)\right)$ arcs.
	
	\medskip
	
	\textbf{Case 3:} \textit{$D$ does not contain a dominating vertex, nor does it contain a dominating $\vec{C_3}$.}
	Since $D\in \vEX(k,n)$, we obtain that
	\[
		D = D_1 \dis D_2 \dis \dots \dis D_s
	\]
	where $\overline{D}_1, \overline{D}_2, \ldots, \overline{D}_s$ are the components of $\overline{D}$. Since $n\leq 2k-2$, it follows from Theorem~\ref{theorem:decomposable} that $s\geq 2$.
	 For every $i\in [1,s]$, we obtain that $D_i\in \vCR(k_i,n_i)$ (by Theorem~\ref{theorem:dirac-join}) and $D_i$ is not decomposable, which implies that $n_i\geq 2k_i-1$ (by Theorem~\ref{theorem:decomposable}). Let $t$ be the number of indices $i$ for which $n_i = 2k_i-1$. Assume first that $t\leq 1$, then we have
	\[
	n = \sum_{i=1}^s n_i \geq \sum_{i=1}^s 2k_i-1 = 2k-1.
	\]
	This is a contradiction since $n\leq 2k-2$. Assume then that $t\geq 2$. By symmetry, we may assume $n_1 = 2k_1 - 1$ and $n_2 = 2k_2 - 1$. Since $D\in \vEX(k,n)$, we obtain as in the first case, that $D_i\in \vEX(k_i,n_i)$ for $i\in [1,s]$. 
Since $D$ does not contain any dominating vertex nor any dominating $\vec{C_3}$, we have $k_i\geq 3$. Hence we have $n_1\geq k_1+2$ and $n_2\geq k_2+2$. Then the induction hypothesis implies that $D_i$ belongs to $\vDG(k_i)$ for $i\in \{1,2\}$.
On the one hand, this implies that:
	\begin{align*}
		|A(D_1 \dis D_2)| &= 2\left( {\binom{2k_1-1}{2}} - ((k_1-1)^2 + 1) \right)\\
						&+ 2\left( {\binom{2k_2-1}{2}} - ((k_2-1)^2 + 1) \right)\\
						& + 2(2k_1-1)(2k_2-1).
	\end{align*}
	
	On the other hand, the minimality of $D$ implies that $D_1 \dis D_2$ belongs to $\vEX(k_1+k_2,2(k_1+k_2)-2)$. A digraph $D'$ in $K_{1} \dis \vDG(k_1+k_2-1)$ belongs to $\vCR(k_1+k_2,2(k_1+k_2)-2)$ and has $2\left( {\binom{2k_1+2k_2-2}{2}} - ((k_1+k_2-2)^2 + 1) \right)$ arcs.
	We finally reach a contradiction to the minimality of $|A(D_1 \dis D_2)|$, since
	\begin{align*}
        &|A(D_1 \dis D_2)| -|A(D')|\\
		&=|A(D_1 \dis D_2)| - 2\left( {\binom{2k_1+2k_2-2}{2}} - ((k_1+k_2-2)^2 + 1) \right)\\
		&= 4k_1k_2 - 4k_1 - 4k_2 +2 \geq 2.
	\end{align*}
	This concludes the proof.
\end{proof}

Theorem~\ref{theorem:Gallai+digraph} applied to bidirected graphs gives Gallai's result (Theorem~\ref{theorem:gallai2}). Gallai's original proof was much longer, since he did not use Dirac's result from 1974. On the other hand, Gallai's result combined with the Kotochka--Yancey bound for $\ex(k,n)$ implies Dirac's result, since the Kostochka--Yancey bound is better that the Dirac bound if $n\geq 2k$.

\section*{Acknowledgement}
The first author has been financially supported by the research grant DIGRAPHS ANR-19-CE48-0013 and by the French government, through the EUR DS4H Investments in the Future project managed by the National Research Agency (ANR) with the reference number ANR-17-EURE-0004.

\end{document}